\documentclass[reqno, 10pt, a4paper]{amsart}

\usepackage[T1]{fontenc}
\usepackage{amssymb}
\usepackage{hyperref}
\usepackage{graphicx}
\usepackage{tikz-cd}
\usepackage{bbm}
\usepackage{booktabs}
\usepackage{mathrsfs}

\begin{document}

\title{On a result of Green and Griffiths}
\author{Qizheng Yin}
\address{Université Pierre et Marie Curie Paris VI, 4 place Jussieu, 75252 Paris cedex 05, France \newline \indent Universiteit van Amsterdam, Korteweg-de Vries Instituut voor Wiskunde, Postbus 94248, 1090 GE Amsterdam, The Netherlands}
\email{qyin@uva.nl}
\date{\today}

\begin{abstract}
We present a simple proof of a result of Green and Griffiths, which states that for the generic curve $C$ of genus $g \geq 4$, the $0$-cycle $K \times K - (2g - 2) K_{\Delta}$ is non-torsion in ${\rm CH}^2(C \times C)$. The proof is elementary and works in all characteristics. 
\end{abstract}
\maketitle

\section*{\bf Introduction}

\noindent Let $C$ be a smooth algebraic curve of genus $g$ over a field $k$. Denote by $K \in {\rm CH}^1(C)$ the (class of a) canonical divisor of $C$. We consider the $0$-cycle
\begin{equation} \label{defZ}
Z := K \times K - (2g - 2) K_{\Delta} \in {\rm CH}^2(C \times C),
\end{equation}
where $K_{\Delta}$ is the divisor $K$ on the diagonal $\Delta \subset C \times C$. This cycle is of degree $0$ and lies in the kernel of the Albanese map.

It is easy to see that $Z = 0$ when the genus $g = 0, 1, 2$. Faber and Pandharipande showed that  also $Z = 0$ when $g = 3$, using the fact that curves of genus $3$ are either hyper-elliptic or plane curves. They asked if $Z$ vanishes in general. In \cite{GG03}, Green and Griffiths answered this with the following result.

\medskip
\noindent {\it Theorem.} --- If $C$ is the generic curve of genus $g \geq 4$, then $Z \neq 0$ in ${\rm CH}^2(C \times C)_{\mathbb{Q}}$.
\medskip

Their proof is Hodge-theoretic. It involves lengthy calculations of a delicate infinitesimal invariant. In this note we give an elementary proof of this result, which also works in positive characteristic. It consists of two separate steps:

\medskip
\noindent {\it Step 1: a problem on the Jacobian.} --- We observe that $Z$ is symmetric, so it naturally lives on the second symmetric power $C^{[2]}$ of $C$. The latter is closely related to the Jacobian $J$ of $C$ via the map $C^{[2]} \to J$ (with respect to a point $x_0 \in C$). We show that $Z$ is the pull-back of an explicit codimension $2$ cycle $W$ on $J$, which is tautological in the sense of Polishchuk \cite{Pol07}. In particular, the $\mathfrak{sl}_2$-action studied in {\it loc.\,cit.}~gives a unified proof that $Z = 0$ for $g = 3$.

\medskip
\noindent {\it Step 2: a degeneration argument.} --- We consider $W$ in the relative setting and it lives on the universal Jacobian $\pi \colon \mathscr{J} \to \mathscr{M}_{g, 1}$. Although Abel-Jacobi trivial fiberwise, the cycle $W$ gives a class ${\rm cl}(W)$ in $H^2(\mathscr{M}_{g, 1}, R^2\pi_*\mathbb{Q})$ (over $\mathbb{C}$, or $H^2\big(\mathscr{M}_{g, 1}, R^2\pi_*\mathbb{Q}_{\ell}(2)\big)$ in general). If $W$ is trivial on the generic fiber, then there should exist an open subset $U \subset \mathscr{M}_{g, 1}$ such that the restriction of ${\rm cl}(W)$ is zero in $H^2(U, R^2\pi_*\mathbb{Q})$. So it remains to show that such a $U$ does not exist. A key lemma by Fakhruddin ({\it cf.}~\cite{Fak96}, Lemma 4.1) reduces this to an argument on the boundary of $\mathscr{M}_{g, 1}$. There we construct explicit families of stable curves and study the cycle class of $W$. It turns out that even the simplest families of `test curves' will suffice for the proof. 

\medskip
\noindent {\it Philosophical note.} --- It is in general very difficult to detect non-trivial cycles in the kernel of the Abel-Jacobi map. Results in this direction are mostly variational, often obtained by calculating infinitesimal invariants on the generic fiber. The invariants are, essentially, Hodge-theoretic objects associated to certain Leray filtrations, and the calculation is usually difficult.

Now since we work with an Abelian scheme $\mathscr{J}$, the classical Leray filtration is actually a multiplicative decomposition ({\it cf.}~\cite{Voi12},~4.3.3). Its compatibility with the Beauville decomposition in ${\rm CH}(\mathscr{J})_{\mathbb{Q}}$ tells us exactly in which cohomology group lies the cycle class. Finally via the degeneration argument we take full advantage of the boundary of $\mathscr{M}_g$ (or $\mathscr{M}_{g, 1}$), which is missing in the infinitesimal approach.

Although elementary, this method can also be used to detect cycles that lie deeper in the conjectural Bloch-Beilinson filtration ({\it cf.}~\cite{Yin13}).

\medskip
\noindent {\it Notation and conventions.} --- We work over a field $k$ of arbitrary characteristic.  Since the main result is a geometric statement ({\it cf.}~note after Theorem~\ref{mainthm}), we assume $k$ to be algebraically closed. From now on, Chow groups are with $\mathbb{Q}$-coefficients. By a `{\it cycle\/}' we mean the rational equivalence class of a cycle. The word `{\it generic\/}' is taken in the schematic sense. Over $\mathbb{C}$ (or any uncountable field), the term `{\it very general\/}' is often used, which means outside a countable union of Zariski-closed proper subsets of the base scheme ({\it cf.}~Corollary~\ref{verygeneral}).

\medskip
\noindent {\it Acknowledgements.} --- Most of the ideas and techniques we use are due to Alexander Polishchuk and Najmuddin Fakhruddin. I am deeply grateful to my advisor Ben Moonen, who is always there to help. Special thanks to Rahul Pandharipande for inviting me to ETH Zürich and for his interest.

\medskip
\section{\bf Connections with the Jacobian}
\medskip

\subsection{\texorpdfstring{\!\!}{}} \label{defq} We briefly review Polishchuk's work on the tautological ring of a Jacobian $J$ ({\it cf.}~\cite{Pol07}). Let $C$ be a smooth curve of genus $g$ over $k$, and let $(J, \theta)$ be the Jacobian of $C$. Recall the Beauville decomposition ${\rm CH}^i(J) = \oplus_{j = i - g}^i {\rm CH}^i_{(j)}(J)$, with
\begin{equation*}
{\rm CH}^i_{(j)}(J) := \big\{\alpha \in {\rm CH}^i(J) : [n]^*(\alpha) = n^{2i - j} \alpha \, \textrm{ for all } n \in \mathbb{Z}\big\}.
\end{equation*}
We identify $J$ with its dual $J^t$ via the canonical principal polarization, and denote by $\mathscr{P}$ the Poincaré line bundle on $J \times J$. There is the Fourier transform $\mathscr{F} \colon {\rm CH}^i_{(j)}(J) \xrightarrow{\sim} {\rm CH}^{g - i + j}_{(j)}(J)$ given by $\alpha \mapsto {\rm pr}_{2, *} \big({\rm pr}_1^*(\alpha) \cdot {\rm ch}(\mathscr{P})\big)$, where ${\rm pr}_1, {\rm pr}_2$ are the two projections ({\it cf.}~\cite{Bea86}).

Choose a base point $x_0 \in C$, and let $\iota \colon C \hookrightarrow J$ be the embedding given by $x \mapsto \mathscr{O}_C(x - x_0)$. Consider the $1$-cycle $[C] := \big[\iota(C)\big]$ and its components $[C]_{(j)} \in {\rm CH}^{g - 1}_{(j)}(J)$. Define
\begin{align*}
p_i & := \mathscr{F}\big([C]_{(i - 1)}\big) \in {\rm CH}^i_{(i - 1)}(J), \textrm{ for } i \geq 1, \\
q_i & := \mathscr{F}\big(\theta \cdot [C]_{(i)}\big) \in {\rm CH}^i_{(i)}(J), \textrm{ for } i \geq 0.
\end{align*}
The $\mathbb{Q}$-subalgebra of ${\rm CH}(J)$ generated by $\{p_i\}_{i \geq 1}$ and $\{q_i\}_{i \geq 0}$ is called the tautological ring of $J$, denoted by $\mathscr{T}(J)$. Polishchuk proved that it is stable under $\mathscr{F}$ and the Pontryagin product `$*$'.

There is a natural $\mathfrak{sl}_2$-action on ${\rm CH}(J)$ and on $\mathscr{T}(J)$ (here $\mathfrak{sl}_2 = \mathbb{Q} \cdot e + \mathbb{Q} \cdot f + \mathbb{Q} \cdot h$), defined by
\begin{equation*}
e(\alpha) := p_1 \cdot \alpha, \quad f(\alpha) := -[C]_{(0)} * \alpha,
\end{equation*}
and
\begin{equation*}
h(\alpha) := (2i - j - g) \alpha, \textrm{ for } \alpha \in {\rm CH}^i_{(j)}(J).
\end{equation*}
Polishchuk also showed that $f$ acts on $\mathscr{T}(J)$ via the differential operator $\mathscr{D}$ given by
\begin{equation} \label{diffop}
\mathscr{D} := -\frac{1}{2}\sum_{i, j \geq 1}\binom{i + j}{j}p_{i + j - 1}\partial_{p_i}\partial_{p_j} - \sum_{i, j \geq 1}\binom{i + j - 1}{j}q_{i + j - 1}\partial_{q_i}\partial_{p_j} + \sum_{i \geq 1}q_{i - 1}\partial_{p_i}.
\end{equation}

Now consider the map $\phi \colon C \times C \to J$ given by $(x, y) \mapsto \mathscr{O}_C(x + y - 2x_0)$. We would like to express the cycle $Z \in {\rm CH}^2(C \times C)$ in (\ref{defZ}) as the pull-back of a certain cycle $W \in \mathscr{T}^2(J)$ under $\phi$. Since $Z$ is Abel-Jacobi trivial, we should look for $W$ in $\mathscr{T}^2_{(2)}(J)$, which is spanned by $q_1^2$ and $q_2$.

\subsection{\it Proposition} --- {\it We have $Z = \phi^*(W)$ with $W := 2\big(q_1^2 - (2g - 2)q_2\big) \in \mathscr{T}^2_{(2)}(J)$.}

\begin{proof}
This is done by an explicit calculation. The essential ingredients are the pull-back of $\theta$ via $\iota \colon C \to J$, and the pull-back of $c_1(\mathscr{P})$ via $(\iota, \phi) \colon C \times (C \times C) \to J \times J$. Write $\xi := \iota^*(\theta)$ and $\ell := (\iota, \phi)^*\big(c_1(\mathscr{P})\big)$. Then we have
\begin{equation*}
\phantom{\textit{cf.}~\textrm{[Pol07], (1,1)}} \xi = \frac{1}{2}K + [x_0] \tag{{\it cf.}~\cite{Pol07}, (1.1)}
\end{equation*}
and
\begin{equation*}
\phantom{\textit{cf.}~\textrm{[Pol07]b, (2,1)}} \ell = [\Delta_{1}] + [\Delta_{2}] - 2{\rm pr}_1^*\big([x_0]\big) - {\rm pr}_2^*\big([x_0 \times C] + [C \times x_0]\big), \tag{{\it cf.}~\cite{Pol07b}, (2.1)}
\end{equation*}
where $\Delta_1 = \big\{(x, x, y) : x, y \in C\big\}$ and $\Delta_2 = \big\{(x, y, x) : x, y \in C\big\}$, and ${\rm pr}_1 \colon C \times (C \times C) \to C$, ${\rm pr}_2 \colon C \times (C \times C) \to C \times C$ are the projections.

By chasing through several cartesian squares, we have
\begin{align*}
&~\phi^*\Big(\mathscr{F}\big(\theta \cdot [C]\big)\Big) = {\rm pr}_{2, *}\big({\rm pr}_1^*(\xi) \cdot \exp(\ell)\big) \\
= &~{\rm pr}_{2, *}\bigg({\rm pr}_1^*\big(\frac{1}{2}K + [x_0]\big) \cdot \exp\Big([\Delta_1] + [\Delta_2] - 2{\rm pr}_1^*\big([x_0]\big)\Big)\bigg) \cdot \exp\big(\!-\![x_0 \times C] - [C \times x_0]\big) \\
= &~{\rm pr}_{2, *}\bigg({\rm pr}_1^*\Big(\big(\frac{1}{2}K + [x_0]\big) \cdot \exp\big(\!-\!2[x_0]\big)\Big) \cdot \exp\big([\Delta_1] + [\Delta_2]\big)\bigg) \cdot \exp\big(\!-\![x_0 \times C] - [C \times x_0]\big) \\
= &~{\rm pr}_{2, *}\Big({\rm pr}_1^*\big(\frac{1}{2}K + [x_0]\big) \cdot \exp\big([\Delta_1] + [\Delta_2]\big)\Big) \cdot \exp\big(\!-\![x_0 \times C] - [C \times x_0]\big).
\end{align*}
Then by expanding the exponentials while keeping track of the codimension, we get
\begin{align*}
\phi^*(q_1) & = {\rm pr}_{2, *}\Big({\rm pr}_1^*\big(\frac{1}{2}K + [x_0]\big) \cdot \big([\Delta_1] + [\Delta_2]\big)\Big) - {\rm pr}_{2, *}{\rm pr}_1^*\big(\frac{1}{2}K + [x_0]\big) \cdot \big([x_0 \times C] + [C \times x_0]\big) \\
& = \frac{1}{2}\big(K \times [C] + [C] \times K\big) - (g - 1)\big([x_0 \times C] + [C \times x_0]\big) \\
\intertext{and}
\phi^*(q_2) & = {\rm pr}_{2, *}\Big({\rm pr}_1^*\big(\frac{1}{2}K + [x_0]\big) \cdot \frac{1}{2}\big([\Delta_1] + [\Delta_2]\big)^2\Big) \\
& \qquad - {\rm pr}_{2, *}\Big({\rm pr}_1^*\big(\frac{1}{2}K + [x_0]\big) \cdot \big([\Delta_1] + [\Delta_2]\big)\Big) \cdot \big([x_0 \times C] + [C \times x_0]\big) \\
& \qquad + {\rm pr}_{2, *}{\rm pr}_1^*\big(\frac{1}{2}K + [x_0]\big) \cdot \frac{1}{2}\big([x_0 \times C] + [C \times x_0]\big)^2 \\
& = {\rm pr}_{2, *}\Big({\rm pr}_1^*\big(\frac{1}{2}K + [x_0]\big) \cdot \big([\Delta_1] \cdot [\Delta_2]\big)\Big) - \frac{1}{2}\big(K \times [x_0] + [x_0] \times K\big) + (g - 2)[x_0 \times x_0] \\
& = \frac{1}{2} K_{\Delta} - \frac{1}{2}\big(K \times [x_0] + [x_0] \times K\big) + (g - 1)[x_0 \times x_0].
\end{align*}
Hence
\begin{align*}
\phi^*(q_1^2) = \frac{1}{2}K \times K - (g - 1)\big(K \times [x_0] + [x_0] \times K\big) + 2(g - 1)^2[x_0 \times x_0],
\end{align*}
and we obtain $\phi^*\Big(2\big(q_1^2 - (2g - 2)q_2\big)\Big) = K \times K - (2g - 2)K_{\Delta}$.
\end{proof}

\subsection{\it Corollary} \label{motive} --- (i) {\it We have $Z = 0$ if and only if $W = 0$. In particular, that $W$ vanishes or not is independent of the point $x_0 \in C$.}

\smallskip
\noindent\makebox[\leftmargin][r]{(ii) }{\it If $g = 3$, then $Z = 0$ in 
${\rm CH}^2(C \times C)$.}

\begin{proof}
For (i), we know that $W \in {\rm CH}\big(h^2(J)\big)$, where $h^2(J)$ is the second component in the motivic decomposition of $h(J)$ ({\it cf.}~\cite{DM91}). The map $\phi$ induces
\begin{equation*}
\phi^* \colon h^2(J) \xrightarrow{\sim} S^2\big(h^1(C)\big) \hookrightarrow h^2(C^{[2]}) \hookrightarrow h^2(C \times C),
\end{equation*}
which realizes $h^2(J)$ as a direct summand of $h^2(C \times C)$. Therefore $\phi^*|_{{\rm CH}(h^2(J))}$ is injective.

For (ii), consider the cycle $p_2^2 \in \mathscr{T}^4_{(2)}(J)$. When $g = 3$, we have $p_2^2 = 0$ for dimension reasons. Apply twice the differential operator $\mathscr{D}$ in~(\ref{diffop}), and we obtain
\begin{equation*}
\mathscr{D}^2(p_2^2) = \mathscr{D}(-6p_3 + 2q_1p_1) = 2(q_1^2 - 4q_2) = 0.
\end{equation*}
So $W = 0$, and thus $Z = 0$.
\end{proof}

\subsection{\it Remarks} --- (i) By a classical theorem of Ro\u{\i}tman \cite{Roi80}, the vanishing of $Z$ with $\mathbb{Q}$-coefficients implies its vanishing with $\mathbb{Z}$-coefficients.

\smallskip
\noindent\makebox[\leftmargin][r]{(ii) }It would be interesting to study the vanishing locus of $Z$ when $g \geq 4$. Also conjecturally $Z$ vanishes if the curve $C$ is defined over $\overline{\mathbb{Q}}$.

\medskip
\section{\bf Fakhruddin's degeneration argument}
\medskip

\subsection{\texorpdfstring{\!\!}{}} We switch to the relative setting and we include stable curves of compact type. Let $S$ be a smooth connected variety of dimension $d$ over $k$. Let $p \colon \mathscr{C} \to S$ be a family of stable $1$-pointed curves of genus $g$ and of compact type. Denote by $\pi \colon \mathscr{J} \to S$ the relative Jacobian of $\mathscr{C}$, which is an Abelian scheme. We write $x_0 \colon S \to \mathscr{C}$ for the section (marked point) of $\mathscr{C}$, and $\sigma_0 \colon S \to \mathscr{J}$ for the zero section of $\mathscr{J}$. As in the absolute case, there is an embedding $\iota \colon \mathscr{C} \hookrightarrow \mathscr{J}$ that gives us a diagram.
\begin{equation} \label{diagram}
\begin{tikzcd}[row sep=normal, column sep=small]
\mathscr{C} \arrow[hookrightarrow]{rr}{\iota} \arrow{dr}{p} & & \mathscr{J} \arrow{dl}[swap]{\pi} \\
& S \arrow[bend left=50]{ul}{x_0} \arrow[bend right=50]{ur}[swap]{\sigma_0}
\end{tikzcd}
\end{equation}

\subsection{\texorpdfstring{\!\!}{}} On ${\rm CH}(\mathscr{J})$, we again have a decomposition ${\rm CH}^i(\mathscr{J}) = \oplus{\rm CH}^i_{(j)}(\mathscr{J})$ such that $[n]^*$ is the multiplication by $n^{2i - j}$ on ${\rm CH}^i_{(j)}(\mathscr{J})$ ({\it cf.}~\cite{DM91}). Regarding the $\ell$-adic realization of the motive $h(\mathscr{J}/S)$, there is a canonical decomposition
\begin{equation*}
\phantom{\textrm{(`$(r)$' for Tate twists)}} R\pi_*\mathbb{Q}_{\ell}(r) \simeq \sum_i R^i\pi_*\mathbb{Q}_{\ell}(r)[-i], \tag{`$(r)$' for Tate twists}
\end{equation*}
with $[n]^*$ acting on $R^i\pi_*\mathbb{Q}_{\ell}(r)$ by the multiplication by $n^i$ ({\it cf.}~\cite{Del68}, 2.19). This decomposition is compatible with the multiplicative structure
\begin{equation*}
R\pi_*\mathbb{Q}_{\ell}(r) \otimes R\pi_*\mathbb{Q}_{\ell}(r') \to R\pi_*\mathbb{Q}_{\ell}(r + r')
\end{equation*}
given by the cup product ({\it cf.}~\cite{Voi12},~4.3.3). It follows that we have a multiplicative decomposition
\begin{equation*}
H^m\big(\mathscr{J}, \mathbb{Q}_{\ell}(r)\big) \simeq \bigoplus_{i + j = m} H^j\big(S, R^i\pi_*\mathbb{Q}_{\ell}(r)\big).
\end{equation*}
Comparing the action of $[n]$ on Chow groups and on cohomology, we know that the cycle class map ${\rm cl} \colon {\rm CH}^i(\mathscr{J}) \to H^{2i}\big(\mathscr{J}, \mathbb{Q}_{\ell}(i)\big)$ decomposes as a sum of maps
\begin{equation} \label{cycleclass}
{\rm cl} \colon {\rm CH}^i_{(j)}(\mathscr{J}) \to H^j\big(S, R^{2i - j}\pi_*\mathbb{Q}_{\ell}(i)\big),
\end{equation}
which respect the multiplicative structures on both sides. Note that if the base field $k = \mathbb{C}$, one may work with singular cohomology with coefficients in $\mathbb{Q}$.

\subsection{\texorpdfstring{\!\!}{}} Consider a cycle $\alpha \in {\rm CH}^i_{(j)}(\mathscr{J})$. Denote by $J_{\eta}$ the generic fiber of $\mathscr{J} \to S$, and by $\alpha_{\eta} \in {\rm CH}^i_{(j)}(J_{\eta})$ the restriction of $\alpha$ to $J_{\eta}$. Now suppose $\alpha_{\eta} = 0$.  By the `spreading-out' procedure ({\it cf.}~\cite{Voi12},~2.1), there exists a non-empty open subset $U \subset S$ such that $\alpha_U = 0 \in {\rm CH}^i_{(j)}(\mathscr{J}_U)$, where $\mathscr{J}_U := \mathscr{J} \times_S U$ and $\alpha_U := \alpha|_{\mathscr{J}_U}$. Combining with the cycle class map~(\ref{cycleclass}), we have the following implication.

\subsection{\it Proposition} --- \label{spread} {\it If $\alpha_{\eta} = 0$, then there exists a non-empty open subset $U \subset S$ such that}
\begin{equation*}
{\rm cl}(\alpha_U) = 0 \in H^j\big(U, R^{2i - j}\pi_*\mathbb{Q}_{\ell}(i)\big). \tag*{\qed}
\end{equation*}

\smallskip
We consider the cycles $q_1^2, q_2$ and $W$ ({\it cf.}~Section~\ref{defq}) in the relative setting (\ref{diagram}). More precisely, denote by $\theta \in {\rm CH}^1_{(0)}(\mathscr{J})$ the divisor class
corresponding to the canonical principal polarization $\lambda \colon \mathscr{J} \xrightarrow{\sim} \mathscr{J}^t$ (so $2\theta$ is the pull-back of the first Chern class of the Poincaré bundle $\mathscr{P}$ under the map ${\rm id} \times \lambda \colon \mathscr{J} \to \mathscr{J} \times \mathscr{J}^t$). Again we identify $\mathscr{J}$ with $\mathscr{J}^t$ and regard the Fourier transform $\mathscr{F}$ as an endomorphism of ${\rm CH}(\mathscr{J})$. Generalizing the definitions of Section~\ref{defq}, we write $[\mathscr{C}] := \big[\iota(\mathscr{C})\big]$ and let
\begin{equation*}
q_i : = \mathscr{F}\big(\theta \cdot [\mathscr{C}]_{(i)}\big) \in {\rm CH}^i_{(i)}(\mathscr{J}), \textrm{ for } i \geq 0.
\end{equation*}
If the context requires it, we shall denote these cycles by $q_i(\mathscr{C})$. As before, we define $W := 2\big(q_1^2 - (2g - 2)q_2\big) \in {\rm CH}^2_{(2)}(\mathscr{J})$. 

Our main focus is the case $S = \mathscr{M}_{g, 1}^{\rm ct}$, {\it i.e.}~the moduli stack of stable $1$-pointed curves of genus $g$ and of compact type. The fact that $\mathscr{M}_{g,1}^{\rm ct}$ is a stack plays no role in the discussion. In fact, since the Chow groups are with $\mathbb{Q}$-coefficients, for our purpose ({\it cf.}~Theorem~\ref{mainthm}) it is equivalent to work over a finite cover of the moduli stack that is an honest variety.

The goal is to prove that for $g \geq 4$, we have $W \neq 0$ generically over $\mathscr{M}_{g, 1}^{\rm ct}$. After Proposition~\ref{spread}, we would like to show that for all non-empty open subsets $U \subset \mathscr{M}_{g, 1}^{\rm ct}$, we have
\begin{equation*}
{\rm cl}(W_{U}) \neq 0 \in H^2\big(U, R^2\pi_*\mathbb{Q}_{\ell}(2)\big).
\end{equation*}
Using the following lemma by Fakhruddin we can reduce the proof of this to a calculation on the boundary of $\mathscr{M}_{g, 1}^{\rm ct}$.

\subsection{\it Lemma {\rm (\texorpdfstring{\cite{Fak96}}{[Fak96]}, Lemma~4.1)}} --- {\it Let $X, S$ be smooth connected varieties over $k$ and $\pi \colon X \to S$ be a smooth proper map. Consider a class $h \in H^m\big(X, \mathbb{Q}_{\ell}(r)\big)$. Suppose there exists a non-empty subvariety $T \subset S$ such that for all non-empty open subsets $V \subset T$, we have $h_V \neq 0$, where $h_V := h|_{X_V}$. Then for all non-empty open subsets $U \subset S$, we have $h_U \neq 0$. \qed}

\medskip
Therefore to achieve the goal, it suffices to construct a family of `test curves' over a variety $T$ on the boundary of $\mathscr{M}_{g, 1}^{\rm ct}$, and to show that the class of $W$ does not vanish over any non-empty open subset of $T$. In fact, we can prove a slightly stronger result.

\subsection{\it Theorem} \label{mainthm} --- {\it When $g \geq 4$, the cycles $q_1^2$ and $q_2$ are linearly independent on the generic Jacobian (over $\mathscr{M}_{g, 1}^{\rm ct}$). In particular, we have $W \neq 0$ on the generic Jacobian.}

\medskip
Note that Theorem~\ref{mainthm} is of geometric nature: if the statement is true over the base field $k$, then it is automatically true over any base field $k' \subset k$. Therefore the theorem still holds over an arbitrary field (not necessarily algebraically closed). Together with Corollary~\ref{motive}~(i), it implies the result of Green and Griffiths.

The rest of this paper is devoted to the construction of the `test curves' and the proof of Theorem~\ref{mainthm}. We shall construct two families of curves over the same base scheme $T$. We will show that for any non-trivial linear combination of $q_1^2$ and $q_2$, at least one of the two families will give a cohomology class that does not vanish when restricted to non-empty open subsets of $T$. For simplicity, we begin with the case $g = 4$, while the proof for the general case is almost identical ({\it cf.}~Section~\ref{generalcase}).

\subsection{\it Case \texorpdfstring{$g = 4$}{g = 4}} --- Take two smooth curves $C_1$ and $C_2$ of genus $2$ over $k$, with Jacobians $(J_1, \theta_1)$ and $(J_2, \theta_2)$. Let $x$ (resp.~$y$) be a varying point on $C_1$ (resp.~$C_2$), and $c$ be a fixed point on $C_2$. We construct the first family of stable curves by joining $x$ and $y$ and using $c$ as the marked point, and then the second family by joining $x$ and $c$ and using $y$ as the marked point, as is shown in the picture below.
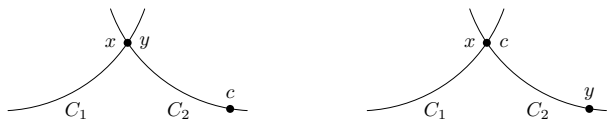
\begin{center}
\includegraphics[height=.065\textheight]{GGPics_1.mps} \quad\quad\quad\quad
\includegraphics[height=.065\textheight]{GGPics_2.mps}
\end{center}

With $x$ and $y$ varying, both families have the same base scheme $T := C_1 \times \big(C_2 \backslash \{c\}\big)$. We denote them by $\mathscr{C} \to T$ and $\mathscr{C}' \to T$ respectively. Observe that $\mathscr{C}$ and $\mathscr{C}'$ have also the same relative Jacobian $\mathscr{J} := J_1 \times J_2 \times T$, a constant Abelian scheme over $T$ via the last projection.

\subsection{\texorpdfstring{\!\!}{}} \label{embed} Consider the embeddings $\mathscr{C} \hookrightarrow \mathscr{J}$ with respect to $c$, and $\mathscr{C}' \hookrightarrow \mathscr{J}$ with respect to $y$. An important fact is that both embeddings naturally extend over $C_1 \times C_2 \supset T$. More precisely, we have
\begin{align*}
\psi_1 \colon C_1 \times C_1 \times C_2 \hookrightarrow J_1 \times J_2 \times C_1 \times C_2 & \, \textrm{ given by } \, (z, x, y) \mapsto \big(\mathscr{O}_{C_1}(z - x), \mathscr{O}_{C_2}(y - c), x, y\big), \\
\psi_2 \colon C_2 \times C_1 \times C_2 \hookrightarrow J_1 \times J_2 \times C_1 \times C_2 & \,\textrm{ given by } \, (w, x, y) \mapsto \big(0, \mathscr{O}_{C_2}(w - c), x, y\big), \\
\psi_1' \colon C_1 \times C_1 \times C_2 \hookrightarrow J_1 \times J_2 \times C_1 \times C_2 & \, \textrm{ given by } \, (z, x, y) \mapsto \big(\mathscr{O}_{C_1}(z - x), \mathscr{O}_{C_2}(c - y), x, y\big), \\
\psi_2' \colon C_2 \times C_1 \times C_2 \hookrightarrow J_1 \times J_2 \times C_1 \times C_2 & \,\textrm{ given by } \, (w, x, y) \mapsto \big(0, \mathscr{O}_{C_2}(w - y), x, y\big).
\end{align*}
We take $\overline{T} := C_1 \times C_2$ as the base scheme and view the other schemes as $\overline{T}$-schemes through the projections onto the last two factors. We also write $\overline{\mathscr{J}} := J_1 \times J_2 \times \overline{T}$. 

Let $\overline{\mathscr{C}} \subset \overline{\mathscr{J}}$ be the union of the images of $\psi_1$ and $\psi_2$; similarly, let $\overline{\mathscr{C}'} \subset \overline{\mathscr{J}}$ be the union of the images of $\psi_1'$ and $\psi_2'$. We see that the restriction of $\overline{\mathscr{C}}$ (resp.~$\overline{\mathscr{C}'}$) to $T$ is exactly $\mathscr{C}$ (resp.~$\mathscr{C}'$). Write $\theta := \theta_1 \times [J_2] \times [\overline{T}] + [J_1] \times \theta_2 \times [\overline{T}]$, and we have $\theta \in {\rm CH}^1_{(0)}(\overline{\mathscr{J}})$. Define
\begin{equation*}
\overline{q_i} := \mathscr{F}\big(\theta \cdot [\overline{\mathscr{C}}]_{(i)}\big) \in {\rm CH}^i_{(i)}(\overline{\mathscr{J}}), \quad \overline{q_i}' := \mathscr{F}\big(\theta \cdot [\overline{\mathscr{C}'}]_{(i)}\big) \in {\rm CH}^i_{(i)}(\overline{\mathscr{J}}).
\end{equation*}
Again, the restriction of $\overline{q_i}$ (resp.~$\overline{q_i}'$) to $T$ is exactly $q_i(\mathscr{C})$ (resp.~$q_i(\mathscr{C}')$).

\subsection{\texorpdfstring{\!\!}{}} As $\overline{\mathscr{J}}$ is a constant Abelian scheme over $\overline{T}$, we have a Künneth decomposition
\begin{equation*}
H^m(\overline{\mathscr{J}}) = \bigoplus_{a_1 + b_1 + a_2 + b_2 = m}H^{a_1}(J_1) \otimes H^{b_1}(C_1) \otimes H^{a_2}(J_2) \otimes H^{b_2}(C_2).
\end{equation*}
Here, and in what follows, we omit the coefficients of the cohomology groups. Also on the right hand side we have sorted the factors in the order $J_1$-$C_1$-$J_2$-$C_2$, as this turns out to be convenient in our calculations. Given a class $h \in H^m(\overline{\mathscr{J}})$, we denote by $h^{[a_1, b_1, a_2, b_2]}$ its Künneth component in the indicated degrees.

In this case, the cycle class map (\ref{cycleclass}) takes the form 
\begin{equation} \label{cycleclass2}
{\rm cl} \colon {\rm CH}^i_{(j)}(\overline{\mathscr{J}}) \to \bigoplus_{\substack{a_1 + a_2 = 2i - j \\ b_1 + b_2 = j}}H^{a_1}(J_1) \otimes H^{b_1}(C_1) \otimes H^{a_2}(J_2) \otimes H^{b_2}(C_2).
\end{equation}
So for $\alpha = \overline{q_1}^2, \overline{q_2}, \overline{q_1}'^2, \overline{q_2}' \in {\rm CH}^2_{(2)}(\overline{\mathscr{J}})$, that ${\rm cl}(\alpha)^{[a_1, b_1, a_2, b_2]} \neq 0$ is possible only if $a_1 + a_2 = 2$ and $b_1 + b_2 = 2$. Moreover, we remark that $H^2(C_1)$ (resp.~$H^2(C_2)$) is supported on a point of $C_1$ (resp.~$C_2$). As we should like to have the cycle classes after restriction to open subsets $V \subset T \subset \overline{T}$, the only interesting components are ${\rm cl}(\alpha)^{[a_1, 1, a_2, 1]}$ with $a_1 + a_2 = 2$ (in fact, we will see in the proof of Proposition~\ref{calcul} that for ${\rm cl}(\alpha)^{[a_1, 1, a_2, 1]}$ to be non-zero, we also have $a_1 = a_2 = 1$).

The following elementary calculation is the key point in the proof of Theorem~\ref{mainthm}.

\subsection{\it Proposition} \label{calcul}--- {\it There exist non-zero classes}
\begin{align*}
h_1 & \in H^1(J_1) \otimes H^1(C_1) \otimes H^1(J_2) \otimes H^1(C_2) \\
h_2, h_4 & \in H^0(J_1) \otimes H^0(C_1) \otimes H^1(J_2) \otimes H^1(C_2) \\
h_3 & \in H^1(J_1) \otimes H^1(C_1) \otimes H^0(J_2) \otimes H^0(C_2)
\end{align*}
{\it such that}
\begin{equation*}
\begin{aligned}
{\rm cl}(\overline{q_2})^{[1, 1, 1, 1]} & = h_1, \\
{\rm cl}(\overline{q_2}')^{[1, 1, 1, 1]} & = -h_1,
\end{aligned}
\quad\quad
\begin{aligned}
{\rm cl}(\overline{q_1}^2)^{[1, 1, 1, 1]} & = 2h_2 \!\smallsmile\! h_3, \\
{\rm cl}(\overline{q_1}'^2)^{[1, 1, 1, 1]} & = -2h_2 \!\smallsmile\! h_3 + 2h_3 \!\smallsmile\! h_4.
\end{aligned}
\end{equation*}
{\it Moreover, the classes $h_2 \!\smallsmile\! h_3$ and $h_3 \!\smallsmile\! h_4$ are also non-zero.}

\begin{proof}
The proof is just a careful analysis of the embeddings $\psi_1, \psi_2, \psi_1', \psi_2'$ defined in Section~\ref{embed}. We first calculate the relevant Künneth components of ${\rm cl}\big([\overline{\mathscr{C}}]_{(i)}\big)$ and ${\rm cl}\big([\overline{\mathscr{C}'}]_{(i)}\big)$. Then by intersecting with ${\rm cl}(\theta)$ and applying $\mathscr{F}$ in cohomology, we obtain the relevant components of ${\rm cl}(\overline{q_i})$ and ${\rm cl}(\overline{q_i}')$.

We start with the cycle classes of $[\overline{\mathscr{C}}]_{(1)}$ and $[\overline{\mathscr{C}}]_{(2)}$. Observe that the image of $\psi_2$ only gives a class in $H^4(J_1) \otimes H^0(C_1) \otimes H^2(J_2) \otimes H^0(C_2)$, which by (\ref{cycleclass2}), does not contribute to either $[\overline{\mathscr{C}}]_{(1)}$ or $[\overline{\mathscr{C}}]_{(2)}$. Regarding $\psi_1$, we may view it as the product of 
\begin{equation*}
\begin{aligned}
\psi_3 \colon C_1 \times C_1 & \hookrightarrow J_1 \times C_1 \\
(z, x) & \mapsto \big(\mathscr{O}_{C_1}(z - x), x\big),
\end{aligned}
\quad\quad
\begin{aligned}
\phantom{C_2 \times} \psi_4 \colon C_2 & \hookrightarrow J_2 \times C_2 \\
y & \mapsto \big(\mathscr{O}_{C_2}(y - c), y\big).
\end{aligned}
\end{equation*}

The class of ${\rm Im}(\psi_3)$ has components in $H^2(J_1) \otimes H^0(C_1)$, $H^1(J_1) \otimes H^1(C_1)$ and $H^0(J_1) \otimes H^2(C_1)$. The third component is irrelevant due to the appearance of $H^2(C_1)$. We claim that the other two components are both non-zero. To see the first, we regard $J_1 \times C_1$ as a constant family over $C_1$. Then $C_1 \times C_1$ is fiberwise an ample divisor, which gives a non-zero class in $H^2(J_1) \otimes H^0(C_1)$. For the component in $H^1(J_1) \otimes H^1(C_1)$, we consider 
\begin{center}
{\renewcommand{\arraystretch}{1.2}\renewcommand{\tabcolsep}{2pt}
\begin{tabular}{@{} c c c c c c c @{}}
$C_1 \times C_1$ & $\xrightarrow{({\rm id}, \Delta)}$ & $C_1 \times C_1 \times C_1$ & $\xrightarrow{(\sigma, {\rm id})}$ & $C_1^{[2]} \times C_1$ & $\xrightarrow{(\varphi, {\rm id})}$ & $J_1 \times C_1$ \\
$(z, x)$ & $\mapsto$ & $(z, x, x)$ & $\mapsto$ & $\big((z, x), x\big)$ & $\mapsto$ & $\big(\mathscr{O}_{C_1}(z + x - 2x), x\big).$
\end{tabular}}
\end{center}
The class of the diagonal in $C_1 \times C_1$ has a component in $H^1(C_1) \otimes H^1(C_1)$ which, viewed as a correspondence, gives the identity $H^1(C_1) \xrightarrow{\sim} H^1(C_1)$. It follows that the class of ${\rm Im}({\rm id}, \Delta)$ has a non-zero component in $H^0(C_1) \otimes H^1(C_1) \otimes H^1(C_1)$. Moreover, we have isomorphisms 
\begin{equation*}
\sigma_* \colon H^0(C_1) \otimes H^1(C_1) \xrightarrow{\sim} H^1(C_1^{[2]}), \quad \varphi_* \colon H^1(C_1^{[2]}) \xrightarrow{\sim} H^1(J_1),
\end{equation*}
the latter due to the fact that $C_1^{[2]}$ is obtained by blowing up a point in $J_1$. Therefore ${\rm Im}(\psi_3)$ as a correspondence gives an isomorphism $H^1(J_1) \xrightarrow{\sim} H^1(C_1)$, which implies a non-zero component in $H^1(J_1) \otimes H^1(C_1)$.

Similarly, the class of ${\rm Im}(\psi_4)$ has non-zero components in $H^4(J_2) \otimes H^0(C_2)$ and $H^3(J_2) \otimes H^1(C_2)$. Now we collect all non-zero contributions to the classes of $[\overline{\mathscr{C}}]_{(1)}$ and $[\overline{\mathscr{C}}]_{(2)}$ that do not involve either $H^2(C_1)$ or $H^2(C_2)$. For $[\overline{\mathscr{C}}]_{(2)}$, there is only one non-zero class
\begin{equation*}
h_1^0 \in H^1(J_1) \otimes H^1(C_1) \otimes H^3(J_2) \otimes H^1(C_2).
\end{equation*}
By intersecting with ${\rm cl}(\theta)$ and applying $\mathscr{F}$, we obtain a non-zero class
\begin{equation*}
h_1 \in H^1(J_1) \otimes H^1(C_1) \otimes H^1(J_2) \otimes H^1(C_2),
\end{equation*}
For $[\overline{\mathscr{C}}]_{(1)}$, there are two non-zero classes
\begin{equation*}
h_2^0 \in H^2(J_1) \otimes H^0(C_1) \otimes H^3(J_2) \otimes H^1(C_2), \quad h_3^0 \in H^1(J_1) \otimes H^1(C_1) \otimes H^4(J_2) \otimes H^0(C_2).
\end{equation*}
Again by intersecting with ${\rm cl}(\theta)$ and applying $\mathscr{F}$, we obtain non-zero classes
\begin{equation*}
h_2 \in H^0(J_1) \otimes H^0(C_1) \otimes H^1(J_2) \otimes H^1(C_2), \quad h_3 \in H^1(J_1) \otimes H^1(C_1) \otimes H^0(J_2) \otimes H^0(C_2).
\end{equation*}
It follows that ${\rm cl}(\overline{q_2})^{[1, 1, 1, 1]} = h_1$ and ${\rm cl}(\overline{q_1}^2)^{[1, 1, 1, 1]} = h_2 \!\smallsmile\! h_3 + h_3 \!\smallsmile\! h_2 = 2h_2 \!\smallsmile\! h_3$.

For the cohomology classes of $\overline{q_2}'$ and $\overline{q_1}'^2$, we remark that the embedding $\psi_1'$ differs from $\psi_1$ only by an action of $[-1]$ on the $J_2$ factor. As a consequence, by repeating the same procedure we obtain classes $h_1' = -h_1$, $h_2' = -h_2$ and $h_3' = h_3$, so that $2h_2' \!\smallsmile\! h_3' = -2h_2 \!\smallsmile\! h_3$. However, this time the embedding $\psi_2'$ makes an additional contribution. The class of ${\rm Im}(\psi_2')$ has a non-zero component
\begin{equation*}
h_4^0 \in H^4(J_1) \otimes H^0(C_1) \otimes H^1(J_2) \otimes H^1(C_2),
\end{equation*}
which belongs to the class of $[\overline{\mathscr{C}'}]_{(1)}$. By intersecting with ${\rm cl}(\theta)$ and applying $\mathscr{F}$, we get a non-zero~class
\begin{equation*}
h_4 \in H^0(J_1) \otimes H^0(C_1) \otimes H^1(J_2) \otimes H^1(C_2).
\end{equation*}
Therefore we have ${\rm cl}(\overline{q_2}')^{[1, 1, 1, 1]} = -h_1$ and ${\rm cl}(\overline{q_1}'^2)^{[1, 1, 1, 1]} = -2h_2 \!\smallsmile\! h_3 + 2h_3 \!\smallsmile\! h_4$.

Finally, since the $0$-th cohomology groups $H^0(C_i)$ and $H^0(J_i)$ are generated by the unit of the ring structures, we see that both $h_2 \!\smallsmile\! h_3$ and $h_3 \!\smallsmile\! h_4$ are non-zero.
\end{proof}

\subsection{\texorpdfstring{\!\!}{}} As $h_1 \neq 0$ and $h_3 \!\smallsmile\! h_4 \neq 0$, it follows from Proposition~\ref{calcul} that that for any $(r, s) \neq (0, 0) \in \mathbb{Q}^2$, at least one of ${\rm cl}(r \cdot \overline{q_1}^2 + s \cdot \overline{q_2})^{[1, 1, 1, 1]}$ and ${\rm cl}(r \cdot \overline{q_1}'^2 + s \cdot \overline{q_2}')^{[1, 1, 1, 1]}$ is non-zero in $H^1(C_1) \otimes H^1(J_1) \otimes H^1(C_2) \otimes H^1(J_2)$.

It remains to ensure that this non-zero cohomology class does not vanish when restricted to non-empty open subsets of $\overline{T} = C_1 \times C_2$, {\it i.e.}~that it is not supported on a divisor of $C_1 \times C_2$. We can achieve this by imposing additional assumptions on $C_1$ and $C_2$. In positive characteristic, we choose $C_1$ to be ordinary and $C_2$ supersingular. Over $\overline{\mathbb{Q}}$, and hence for any $k = \overline{k}$ of characteristic $0$, we take $C_1$ and $C_2$ such that $J_1$ and $J_2$ are both simple, and such that ${\rm End}(J_1) = \mathbb{Z}$ and $J_2$ is of CM type ({\it cf.}~\cite{CF96},~Chapters 14,~15 for explicit examples). In both situations we have ${\rm Hom}(J_1, J_2) = 0$, which implies that there is no non-zero divisor class in $H^1(C_1) \otimes H^1(C_2)$. This completes the proof for $g = 4$.

\subsection{\it General case: end of proof} \label{generalcase} --- When $g >4$, we may attach to both families a constant curve $C_0$ of genus $g - 4$ via a fixed point $c' \in C_0$, and use another fixed point $c'' \in C_0$ as the marked point.
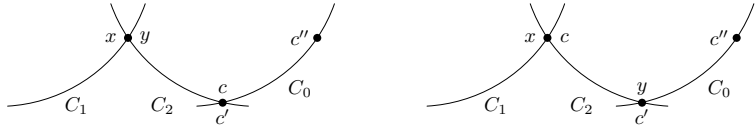
\begin{center}
\includegraphics[height=.07\textheight]{GGPics_3.mps} \quad\quad\quad
\includegraphics[height=.07\textheight]{GGPics_4.mps}
\end{center}
We repeat the same procedure, and the proof is exactly the same. \qed

\subsection{\it Corollary} \label{verygeneral} --- {\it When the base field $k$ is uncountable ({\it e.g.}~$k = \mathbb{C}$) and when $g \geq 4$, the same statement as in Theorem~\ref{mainthm} holds for the Jacobian of a very general curve (over $\mathscr{M}_{g, 1}$).}

\begin{proof}
This is a consequence of the following fact. Let $X, S$ be smooth connected quasi-projective varieties over $k$ and $\pi \colon X \to S$ be a smooth projective map. Consider a cycle $\alpha \in {\rm CH}(X)$. Via an argument using relative Hilbert schemes, one can show that the locus $\big\{s \in S : \alpha_s = 0 \in {\rm CH}(X_s)\big\}$ is a countable union of Zariski-closed subset of $S$ ({\it cf.}~\cite{Voi12},~2.1). Therefore, if $k$ is uncountable and $\alpha$ is non-zero on the generic fiber, then $\alpha$ is non-zero on a very general fiber.
\end{proof}

\end{document}